\documentclass{lms}

\usepackage{amsmath,amssymb}
\usepackage{url}
\usepackage[utf8]{inputenc}
\usepackage[T1]{fontenc}
\usepackage{epstopdf}
\usepackage{epsfig}

\usepackage{amssymb,latexsym, amsmath, amsxtra}
\usepackage{graphicx}
\usepackage{array}

\newtheorem{lemma}{Lemma.}[section]
\newtheorem{theorem}[lemma]{Theorem}

\newtheorem{proposition}[lemma]{Proposition}

\newtheorem{question}[lemma]{Question}
\newtheorem{problem}[lemma]{Problem}

\newcommand{\C}{\mathbb C}

\newcommand{\R}{\mathbb R}
\newcommand{\Z}{\mathbb Z}

\newcommand{\cj}[1]{\overline{#1}}
\newcommand{\abs}[1]{\left| #1 \right|}
\newcommand{\cbr}[1]{\left\{ #1 \right\}}

\newcommand{\spin}{\mathrm{spin}}
\newcommand{\stan}{\mathrm{stan}}
\newcommand{\adj}{\mathrm{adj}}
\newcommand{\Sp}{{\mathrm Sp}}

\def\({\left(}
\def\){\right)}

\numberwithin{equation}{section}

\title{Evaluating $L$-functions with few known coefficients}

\author{
David W. Farmer and Nathan C. Ryan
}

\date{\today}

\classno{11F66}

\begin{document}



\maketitle

\begin{abstract}
We address the problem of evaluating an $L$-function
when only a small number of its Dirichlet coefficients
are known.  
We use the approximate functional equation in a new way and find that it is
possible to evaluate the $L$-function more precisely than one would
expect from the standard approach.  
The method, however, requires considerably more computational effort
to achieve a given accuracy
than would be needed if more Dirichlet coefficients were available.

\end{abstract}

\section{Introduction}

$L$-functions are central to much of contemporary number theory.  Two
celebrated conjectures, the Riemann Hypothesis and the Birch and
Swinnerton-Dyer Conjecture are about values of $L$-functions and were
discovered as a result of the explicit computation of the Riemann
zeta function and the Hasse-Weil $L$-function associated to an elliptic
curve, respectively.  These $L$-functions are, respectively, of
degree one and degree two and it is interesting to verify analogous
conjectures about special values and zeros for higher degree $L$-functions.
Conjectures such as
B\"ocherer's Conjecture \cite{Bocherer2,RyanTornaria} 
and the Bloch-Kato Conjecture
\cite{BlochKato} are about
the central values of degree four and degree three $L$-functions, and the
Grand Riemann Hypothesis asserts that all nontrivial zeros
of an $L$-function, of any degree, lie along the critical line.

In addition to these conjectures, there are a number of other conjectures for the statistical behavior of $L$-functions,
arising from the interplay between random matrix theory and number theoretic heuristics \cite{KS,CFKRS,CFZ}.
One of the main reasons those conjectures are believable is that
large-scale calculations of the value distribution and the
zeros of $L$-functions yield data that support those
conjectures.

The $L$-functions we consider here are associated to Siegel modular forms.
Our examples will use the first non-lift Siegel modular
form on $\Sp(4,\Z)$.  The form has weight 20 and is usually
denoted~$\Upsilon_{20}$.  
Background information beyond what we mention about Siegel modular forms can be found
in~\cite{Abook,Kbook,Skoruppa}. For this paper the relevant information is that a Siegel
modular form is acted on by Hecke operators $T(n)$,
which have eigenvalues~$\lambda(n)$.
It is the eigenvalues $\lambda(p)$ and $\lambda(p^2)$, for $p$ prime,
which are used to define the $L$-functions associated to the modular form.
For $\Upsilon_{20}$ the eigenvalues $\lambda(p)$ have been computed for $p\le 997$,
and the eigenvalues $\lambda(p^2)$ for $p\le 79$~\cite{KohnenKuss}.
These data are available at~\cite{Skoruppa_site}.

There is an $L$-function
$L(s,F, \rho)$ of degree $n$ for
each $n$-dimensional representation $\rho$ of the dual group of
$\textrm{PGSp}(4)$, namely $\textrm{Sp}(4,\C)$.  Associated to a Siegel modular form $F$ is a
sequence of $L$-functions, of degrees 4, 5, 10, 14, 16,~etc.
The degree 4, 5, and 10 $L$-functions are called, respectively,
the spinor, standard, and adjoint, and are denoted
$L(s,F,\spin)$, $L(s,F,\stan)$, and $L(s,F,\adj)$.
Proposition \ref{prop:3 L-functions}, taken from~\cite{FRS}, summarizes the properties
of those $L$-functions for a weight~$k$
Siegel modular form on~$\Sp(4,\Z)$.

The degree 4 and 5 $L$-functions
were shown by Andrianov~\cite{Aspin} and
B\"ocherer~\cite{Bocherer} to
have an analytic continuation and satisfy a functional equation.
The degree 10 $L$-function was
only recently shown by Pitale, Saha, and Schmidt~\cite{PSS} to have an analytic continuation
and satisfy a functional equation.  Those properties for the $L$-functions of degree~14 and
above are still conjectural.  

\subsection{Evaluating $L$-functions}

We are concerned with numerically evaluating $L$-functions.
The standard tool, which is used in available open-source
computational packages \cite{Rub,Dok}
is the approximate functional equation.  See Proposition~\ref{thm:formula}
for the precise formulation.

There are two main difficulties in evaluating high degree $L$-functions.
The first is that if
the $L$-function
$L(s)$ has degree $d$, evaluating $L(\frac12 +it)$ 
using the approximate functional equation requires $\gg (1+|t|)^{d/2}$
Dirichlet series coefficients.  
Here the implied constant depends on the $L$-function and the desired
precision in the answer.
For example, estimating the implied constant from the calculations in
Section~\ref{sec:optimalappfe}, to find the first 1~million zeros of
$L(s,\Upsilon_{20}, \adj)$, the degree 10 $L$-function associated
to~$\Upsilon_{20}$, would require using the approximate functional
equation with around $10^{30}$
Dirichlet series coefficients.  Note that 1~million zeros is not even
a large sample;
for example it is probably not 
sufficient for testing various conjectures about the lower-order terms
in the distribution of spacings between zeros.

The second difficulty is that current methods are incapable of producing a large
number of Dirichlet coefficients of the standard and adjoint $L$-functions
of a Siegel modular form.
The Fourier coefficients indexed by quadratic forms with discriminant up to 3000000
have been computed for $\Upsilon_{20}$~\cite{KohnenKuss}. 
These Fourier coefficients are used to compute the Hecke eigenvalues.
Examination of formulas on page 387 of \cite{Skoruppa} shows that to find
the eigenvalue $\lambda(n)$ of $T(n)$, for $n=p^2$,  
requires the Fourier coefficients indexed by quadratic forms of discriminant up to~$n^2=p^4$.

It gets worse.  By \eqref{eq-satake} and \eqref{eqn:EPs},
the $p$th Dirichlet coefficient of the standard or adjoint $L$-function
requires both $\lambda(p)$ and $\lambda({p^2})$.  Thus, to determine the
first $n$ Dirichlet coefficients of those $L$-functions requires Fourier
coefficients
of the Siegel modular form
of index up to approximately~$n^4$.  The extensive calculations
in~\cite{KohnenKuss} are not even sufficient to determine the 
83rd
Dirichlet coefficient of the standard or adjoint $L$-functions of~$\Upsilon_{20}$.

Of course, one could determine more Dirichlet coefficients by first
finding more Fourier coefficients of the cusp form.
But the $n^4$ relationship makes
this quite expensive, so with current methods it is not feasible to
determine many more Dirichlet coefficients than currently known.  It
is possible that new methods will be developed to determine the Hecke
eigenvalues without extensive computation.  The real problem will
still remain:  how to compute high-degree $L$-functions without
requiring an enormous number of Dirichlet coefficients.

That brings us to the theme of this paper: how accurately can one compute an
$L$-function given a limited number of Dirichlet coefficients. 
As the above discussion indicates, this is a practical problem and there
are many $L$-functions for which it is not currently possible to determine a
reasonably large number of Dirichlet coefficients.

As we describe in this paper, even without knowing many
Dirichlet coefficients, we were able to evaluate the $L$-functions to
surprisingly high accuracy (surprising to us, anyway).  In fact,
we were not able to establish that there is
an absolute limit to the accuracy one can obtain
from only a limited number of coefficients.
However, our method is computationally expensive -- much more
expensive than evaluating the $L$-function in a straightforward
way if more coefficients were available.
If one could find an efficient way to determine the unknown
parameters in our method, that could make it possible to
quickly evaluate high-degree $L$-functions.  See Section~\ref{sec:speculation}
for a discussion.

In the next section we describe the $L$-functions we consider here,
and in Section~\ref{sec:appfe} we recall the approximate functional equation and
how it is used to compute an $L$-function.  In Section~\ref{sec:optimalappfe} we
state the underlying problem and then describe our experiments
evaluating $L(s,\Upsilon_{20},\stan)$ and $L(s,\Upsilon_{20},\adj)$.
In Section~\ref{sec:lp} we describe a second method that performs a
little better than the method in Section~\ref{sec:optimalappfe} and
show how this second method can be used to approximate unknown
Dirichlet series coefficients.

We thank the referee for suggesting we include the material in
Section~\ref{sec:realitychecks}.

\section{The $L$-functions}
The $L$-functions associated to a Siegel modular form are most
conveniently described as Euler products.  The local factors
of the Euler product can be expressed in terms of the Hecke
eigenvalues $\lambda(p)$ and $\lambda(p^2)$, 
but it is more convenient to express them in terms
of the Satake parameters, $\alpha_{0,p}$, $\alpha_{1,p}$, and $\alpha_{2,p}$,
given by
\begin{align}
p^{2k-3} &= \alpha_0^2\alpha_1\alpha_2\\
A &= \alpha_0^2\alpha_1^2\alpha_2+\alpha_0^2\alpha_1\alpha_2^2+\alpha_0^2\alpha_1 + \alpha_0^2\alpha_2 \notag \\
B &= \alpha_0^2\alpha_1^2\alpha_2^2+\alpha_0^2\alpha_1^2+\alpha_0^2\alpha_2^2+\alpha_0^2, \notag 
\end{align}
where
\begin{eqnarray}\label{eq-satake}
\lambda(p)^2 &=& 4p^{2k-3}+2A+B\\\notag
\lambda(p^2) &=& (2-1/p)p^{2k-3} + A +B.
\end{eqnarray}
We suppress the $p$ on the Satake parameters when clear from context.
See \cite{Ryan} for a discussion of how to solve this polynomial
system of three equation for the three unknowns
$\alpha_{0,p},\alpha_{1,p},\alpha_{2,p}$ using Gr\"obner bases.

We rescale the Satake parameters to have the
so-called ``analytic'' normalization $|\alpha_j|=1$,  $ \alpha_0^2\alpha_1\alpha_2 =1$,
which is possible if we assume the Ramanujan bound on the Hecke eigenvalues.
This
corresponds to a simple change of variables in the $L$-functions, so that all our $L$-functions satisfy
a functional equation in the standard form~$s\leftrightarrow 1-s$.

As an error check for the reader who may wish to extend our calculations,
for $\Upsilon_{20}$ we have
$\lambda(2)=-840960$, $\lambda(4)=248256200704$, and the Satake parameters
at~2 are approximately
\begin{equation}
\begin{split}
\{\alpha_0&,\alpha_1,\alpha_2\}=\\
&\{-0.901413+0.43296 i,0.630904 - 0.775861 i,-0.211226+0.977437 i\}.
\end{split}
\end{equation}
Here and throughout this paper, decimal values are truncations
of the true values.

\begin{proposition}\label{prop:3 L-functions}  Suppose $F\in M_k(\Sp(4,\Z))$ is a Hecke eigenform.  Let $\alpha_{0,p}$, $\alpha_{1,p}$,
$\alpha_{2,p}$ be the Satake parameters of $F$ for the prime~$p$,
where we suppress the dependence on $p$ in the formulas below.  For $\rho\in\{\spin,\stan,\adj\}$ we have the $L$-functions $L(s,F,\rho):=\prod_{p \text{ prime}} Q_p(p^{-s},F,\rho)^{-1}$ where
\begin{align}\label{eqn:EPs}
Q_p(X,F,\spin):=\mathstrut &
{(1-\alpha_0 X)(1-\alpha_0\alpha_1 X)
 (1-\alpha_0\alpha_2 X)(1-\alpha_0\alpha_1\alpha_2 X)},\nonumber\\
 Q_p(X,F,\stan):=\mathstrut &(1-X)(1-\alpha_{1}X)(1-\alpha_{1}^{-1}X)
(1-\alpha_{2}X)(1-\alpha_{2}^{-1}X),\nonumber\\
Q_p(X,F,\adj):=\mathstrut &(1- X)^2(1-\alpha_1 X)(1-\alpha_1^{-1} X)
 (1-\alpha_2 X)(1-\alpha_2^{-1} X)\nonumber\\
 &\phantom{X}(1-\alpha_1\alpha_2 X)(1-\alpha_1^{-1}\alpha_2 X)
  (1-\alpha_1\alpha_2^{-1} X)(1-\alpha_1^{-1}\alpha_2^{-1} X),
 \end{align}
give the $L$-series of, respectively, the spinor, standard,  and adjoint $L$-functions.
These $L$-functions satisfy the functional equations:
\begin{align}\label{eqn:FEs}
\Lambda(s,F,\spin):=\mathstrut &
     \Gamma_\C(s+\tfrac12)\Gamma_\C(s+k-\tfrac32)L(s,F,\spin)\nonumber \\
=\mathstrut & (-1)^k \Lambda(1-s,F,\spin),\nonumber\\
\Lambda(s,F,\stan):=\mathstrut &  \Gamma_\R( s)\Gamma_\C(s+k-2)\Gamma_\C(s+k-1) L(s,F,\stan) \nonumber\\
=\mathstrut & \Lambda(1-s,F,\stan)\nonumber,\\
\Lambda(s,F,\adj):=\mathstrut &
\Gamma_\R(s+1)^2\Gamma_\C(s+1)\nonumber\\
 & \times\Gamma_\C(s+k-2)\Gamma_\C(s+k-1)
\Gamma_\C(s+2k-3)  L(s,F,\adj) \nonumber\\
 =\mathstrut  &\Lambda(1-s,F,\adj).
\end{align}
\end{proposition}


In \eqref{eqn:FEs}, we use the normalized $\Gamma$-functions
\[
\Gamma_\R(s):=\pi^{-s/2}\Gamma(s)\ \ \ \ \text{ and }\ \ \ \ \Gamma_\C(s):=2(2\pi)^{-s}\Gamma(s).
\]
The \emph{degree} of an $L$-function is $r+2c$ where $r$ and $c$ are the
number of $\Gamma_\R$ and $\Gamma_\C$ factors in the functional equation,
respectively.
The spin, standard, and adjoint $L$-functions described above are of
degree 4, 5, and~10. The Ramanujan bound for a degree~$d$
$L$-function with Dirichlet series $\sum_{n\geq 1} b_nn^{-s}$  is given by:
\begin{equation}\label{eqn:ram}
|b_{p^j}|\le
\left(
\genfrac{}{}{0pt}{}{d+j-1}{j}
\right).
\end{equation}
Note that this is equivalent to the assertion that the Satake parameters
satisfy~$|\alpha_{j,p}|\le 1$.

\section{The approximate functional equation}\label{sec:appfe}

In this section we describe the  approximate functional equation,
which is the primary tool used to evaluate $L$-functions.
The approximate functional equation involves a test function which
can be chosen with some freedom.  This will play a key role
in our calculations.

\subsection{Smoothed approximate functional equations}
The material in this section is taken from Section~3.2 of~\cite{Rub}.

Let
\begin{equation}
   L(s) = \sum_{n=1}^{\infty} \frac{b_n}{n^s}
\end{equation}
be a Dirichlet series that converges absolutely in a half plane, $\Re(s) > \sigma_1$.

Let
\begin{equation}
    \label{eqn:lambda}
    \Lambda(s) = Q^s
                 \left( \prod_{j=1}^a \Gamma(\kappa_j s + \lambda_j) \right)
                 L(s),
\end{equation}
with $Q,\kappa_j \in {\mathbb{R}}^+$, $\Re(\lambda_j) \geq 0$,
and assume that:
\begin{enumerate}
    \item  $\Lambda(s)$ has a meromorphic continuation to all of ${\mathbb{C}}$ with
           simple poles at $s_1,\ldots, s_\ell$ and corresponding
           residues $r_1,\ldots, r_\ell$.
    \item $\Lambda(s) = \varepsilon \cj{\Lambda(1-\cj{s})}$ for some
          $\varepsilon \in {\mathbb{C}}$, $|\varepsilon|=1$.
    \item For any $\sigma_2 \leq \sigma_3$, $L(\sigma +i t) = O(\exp{t^A})$ for some $A>0$,
          as $\abs{t} \to \infty$, $\sigma_2 \leq \sigma \leq \sigma_3$, with $A$ and the constant in
          the `Oh' notation depending on $\sigma_2$ and $\sigma_3$. \label{page:condition 3}
\end{enumerate}

Note that~\eqref{eqn:lambda} expresses the functional equation in 
more general terms than~\eqref{eqn:FEs}, but it is a simple matter
to unfold the definition of~$\Gamma_\R$ and~$\Gamma_\C$.

To obtain a smoothed approximate functional equation with desirable
properties, Rubinstein \cite{Rub} introduces an auxiliary function.
Let $g: \C \to \C$ be an entire function that, for fixed $s$, satisfies
\begin{equation}\label{eqn:gbound}
    \abs{\Lambda(z+s) g(z+s) z^{-1}} \to 0
\end{equation}
as $\abs{\Im{z}} \to \infty$, in vertical strips,
$-x_0 \leq \Re{z} \leq x_0$. The smoothed approximate functional
equation has the following form.
\begin{theorem}\label{thm:formula}
For $s \notin \cbr{s_1,\ldots, s_\ell}$, and $L(s)$, $g(s)$ as above,
\begin{equation}\label{eqn:formula}
  \Lambda(s) g(s) =
         \sum_{k=1}^{\ell} \frac{r_k g(s_k)}{s-s_k}
         + Q^s \sum_{n=1}^{\infty} \frac{b_n}{n^s} f_1(s,n) 
         + \varepsilon Q^{1-s} \sum_{n=1}^{\infty} \frac{\cj{b_n}}{n^{1-s}} f_2(1-s,n)
\end{equation}
where
\begin{align}\label{eqn:mellin}
   f_1(s,n) &:= \frac{1}{2\pi i}
                   \int_{\nu - i \infty}^{\nu + i \infty}
                    \prod_{j=1}^a \Gamma(\kappa_j (z+s) + \lambda_j)
                    z^{-1}
                    g(s+z)
                    (Q/n)^z
                    dz \notag \\
  f_2(1-s,n) &:= \frac{1}{2\pi i}
                   \int_{\nu - i \infty}^{\nu + i \infty}
                    \prod_{j=1}^a \Gamma(\kappa_j (z+1-s) + \cj{\lambda_j})
                    z^{-1}
                    g(s-z)
                    (Q/n)^z
                    dz
\end{align}
with $\nu > \max \cbr{0,-\Re(\lambda_1/\kappa_1+s),\ldots,-\Re(\lambda_a/\kappa_a+s)}$.
\end{theorem}

In our examples $L(s)$ continues to an entire function, so the first
sum in \eqref{eqn:formula} does not appear.  For fixed
$Q,\kappa,\lambda,\varepsilon$, and sequence~$b_n$, and $g(s)$ as described
below, the right side of \eqref{eqn:formula}
can be evaluated to high precision. 

A reasonable choice for the weight function is 
\begin{equation}\label{eqn:test}
g(s)=e^{i b s + c s^2},
\end{equation}
which by Stirling's formula satisfies \eqref{eqn:gbound}
if  $c>0$, or if $c=0$ and $|b| < \pi d/4$, where $d$ is the degree of
the $L$-function.
Rubinstein~\cite{Rub} uses such a weight function with $b$ chosen
to balance the size of the terms in the approximate functional equation,
minimizing the loss in precision in the calculation.  
In this paper we exploit the fact that there are many choices of weight function,
and so there are many ways to evaluate the $L$-function.
We combine those calculations to extract as much information
as possible from the known Dirichlet coefficients.  This idea is
described in the next section.

In the computations we carry out below, we find it more convenient to
use the Hardy $Z$-function in our computations instead of the
$L$-function itself.  The function $Z$ associated to an $L$-function
$L$ is defined by the properties: $Z(\frac12 + it)$ is a smooth function
which is real if $t$ is
real, and $|Z(\frac12 + it)| = |L(\frac12 + it)|$.

\section{Exploiting the test function in  the approximate functional equation}\label{sec:optimalappfe}

If we let $g(s)=1$ and $s=\frac12 + 10 i$ in the approximate
functional equation \eqref{eqn:formula} for the standard (degree 5)
$L$-function of $\Upsilon_{20}$, we get
\begin{align}\label{eqn:ex1}
Z(\tfrac12+10 i, \Upsilon_{20},\stan)=\mathstrut & -1835.424 -395.011 \,b_2 + 1012.179 \,b_3 + 1906.603 \, b_4
+\nonumber\\
&\mathstrut + 2226.503 \,b_5 + \cdots + 6.840 \times 10^{-9} \,b_{82} +
\nonumber\\
&\mathstrut  + 5.132\times 10^{-9} \,b_{83}+\cdots+3.205 \times 10^{-16}
\,b_{149} \nonumber\\
&\mathstrut + 2.564 \times 10^{-16} \,b_{150} + \cdots .
\end{align}

If instead we let $g(s)=e^{-3i s/2}$ and keep $s=\frac12 + 10i$ then we have
\begin{align}\label{eqn:ex2}
Z(\tfrac12+10 i, \Upsilon_{20},\stan)=\mathstrut &
1.66549 + 1.39643 \, b_{2} -0.658439 \, b_{3} + 0.726149 \, b_{4} 
+
\nonumber\\
&\mathstrut -0.88227 \, b_{5} +\cdots + 1.532 \times 10^{-8} \, b_{82}
\nonumber\\
&\mathstrut + 1.271 \times 10^{-8} \, b_{83}  +\cdots+ 3.514 \times  10^{-14}
\, b_{149} \nonumber \\
&\mathstrut + 2.309 \times 10^{-14} \, b_{150}
+ \cdots .
\end{align}

Note that neither of the above expressions appears optimal: the first
involves large coefficients, which will lead to a loss of precision.
In the second, the terms do not decrease as quickly, so one must use
more coefficients to achieve a given accuracy.
An observation that we exploit is the fact that the above are just
two among a large
number of expressions for the value of the $L$-function at $\frac12+10i$.

Recall that by \cite{KohnenKuss,Skoruppa_site}  we know the Satake
parameters of $\Upsilon_{20}$ for all
$p\le 79$.
Therefore we recognize two types of terms in the approximate functional
equation, as illustrated in the above examples.   There are terms for which we
know the Dirichlet coefficients, such as $b_5$, $b_{82}$, and $b_{150}$.
And there are terms with an unknown Dirichlet coefficient,
such as $b_{83}$ or $b_{149}$.
Actually, there is a third type of term, such as $b_{166}=b_{2}b_{83}$
which is ``partially unknown''.
We can estimate the unknown terms 
by applying the Ramanujan bound to the
Dirichlet coefficient, and evaluate everything else precisely.
Thus, once we choose a test function, we can evaluate an $L$-function
at a given point as
\begin{equation}
Z(s) = \text{calculated\_value}(s) \pm \text{error\_estimate}(s),
\end{equation}
where both the calculated value and the error estimate are functions of
the test function and the set of known Dirichlet coefficients.
For later use, we write
\begin{equation}\label{eqn:error}
\text{error\_estimate}(s) =
\sum_{n: \ b_n \text{unknown}} \delta_n(g,s)  b_n 
\end{equation}
where the $\delta_n(g,s)$ is the coefficient of $b_n$ in
\eqref{eqn:formula}.  The product of $\delta_n(g,s)$ and the Ramanujan
bound for $b_n$ is an upper bound for
the error contributed to computation by the unknown coefficient $b_n$.
In the calculations described below, we directly evaluate the
contributions from the first 2000 Dirichlet coefficients.
For $n\le 2000$ we use the calculated value of $\delta_n(g,s)$ and the Ramanujan
bound for $b_n$ to estimate the contribution.  This is the main
source of the error term in~\eqref{eqn:error}.

While there are rigorous bounds for the contribution of the tail to
the error (see, e.g., Propositions 3.7 and 3.9 in \cite{Molin}) we do
not make use of them for two reasons.  First, those general bounds are much
larger than what is actually observed in our examples.  For instance, in \eqref{eqn:ex2}
it appears that by the 150th term the contribution is less than $10^{-13}$,
and this is confirmed by further computation (to thousands of terms),
showing a steady decrease at the expected rate.  But the general bounds of~\cite{Molin}
require about 8000 terms before the predicted contribution drops below $10^{-13}$.
Second, we consider our method to be experimental and, as
such, did not emphasize being so careful with the bounds and
instead, relied on observation and intuition.  We think, as
illustrated in examples below, the fact that we were able to
obtain consistent values for our calculations of $L$-functions
is evidence that the results are correct and in principle could be made rigorous.

Using the known $b_n$ and applying the Ramanujan bound
\eqref{eqn:ram} to~\eqref{eqn:ex1} we get
\begin{equation}
Z(\tfrac12+10 i, \Upsilon_{20},\stan) = 3.03930\,70838 \pm 3.12 \times 10^{-8} .
\end{equation}
And for~\eqref{eqn:ex2} we get
\begin{equation}
Z(\tfrac12+10 i, \Upsilon_{20},\stan) = 3.03930\,70808 \pm 7.10 \times 10^{-8} .
\end{equation}

\begin{figure}[htp]
\begin{center}
\scalebox{0.7}[0.7]{\includegraphics{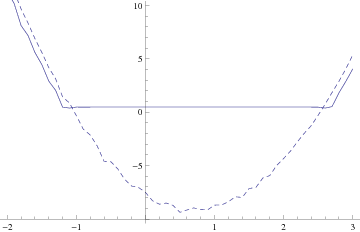}}
\caption{\sf The solid line is the calculated value
and the dashed line is the error estimate in computing
$Z(\frac12+10 i,\Upsilon_{20},\stan)$
 using the available Dirichlet coefficients with the weight function
$g(s)=e^{-i \beta s}$ where $\beta$
is given along the horizontal axis.
The vertical axis is $\log_{10}$ of (the absolute value of) the actual value.
For quite a wide range of test functions, the value of the Z-function
at $\frac12+10i$ is determined with some accuracy, achieving
around 10 decimal digits of accuracy with the optimal choice of $\beta$.
 \label{fig:stan10error}}
\end{center}
\end{figure}

In Figure~\ref{fig:stan10error} we show the
calculated value and error estimate 
for
$Z(\frac12+10 i,\Upsilon_{20},\stan)$ when
evaluated with test
functions of the form~$g(s)=e^{-i\beta s}$.  Note that the vertical
axis in the figure is on a log scale.

As Figure~\ref{fig:stan10error} shows, there is a wide range of
$\beta$ for which it is possible to determine $\Lambda(\frac12+10 i,\Upsilon_{20},\stan)$
with some accuracy. 
With the optimal choice of
$\beta$ the error estimate is approximately $4\times 10^{-10}$.

\begin{figure}[htp] \begin{center}
\scalebox{0.7}[0.7]{\includegraphics{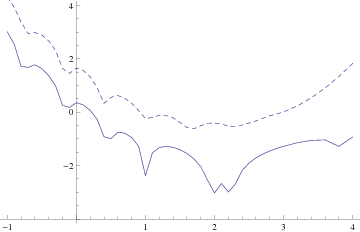}}
\caption{\sf The analogue of Figure~\ref{fig:stan10error}
for $Z(\frac12+5 i,\Upsilon_{20},\adj)$ with the weight function
$g(s)=e^{-i \beta s+\frac{1}{500}(s-5i)^2}$ where $\beta$
 is given along the horizontal axis.
\label{fig:adj5error}
}
\end{center}
\end{figure}

Figure~\ref{fig:adj5error} shows the
calculated value and error estimate
for the adjoint (degree 10) $L$-function
$Z(\frac12+5 i,\Upsilon_{20},\adj)$ when
evaluated with test
functions of the form~$g(s)=e^{-i\beta s+\frac{1}{500}(s-5i)^2}$.

As Figure~\ref{fig:adj5error} shows, every test function of the given
form leads to an error which is larger than the calculated value.
Thus, we can determine that $|Z(\frac12+5 i,\Upsilon_{20},\adj)| < 0.25$,
but with individual test functions of the given form
we cannot even determine if the actual value is positive or negative.

We now introduce a new idea for increasing the accuracy of these
calculations.

\subsection{Optimizing the test function}
In Figure~\ref{fig:stan10error} we see that there are many values of the
parameters
in the test function which give reasonable results.  If there is 
some degree of independence in the errors, then there is hope for obtaining
a smaller error by combining the results of those separate calculations.
Write $Z(s,\Upsilon_{20},\stan,\beta)$ for the output of the
approximate functional equation with weight function~$g(s)=e^{-i \beta s}$.
Consider
\begin{equation}
Z(s,\Upsilon_{20},\stan) = \sum_{j=1}^J c_{\beta_j} Z(s,\Upsilon_{20},\stan,\beta_j)
\end{equation}
where $\sum c_{\beta_j} =1$.  We make the specific choices
\begin{align}\label{eqn:5weights}
s=\mathstrut & \tfrac12+10 i\nonumber\\
(\beta_1,\beta_2,\beta_3,\beta_4,\beta_5) 
   =\mathstrut & \left(\tfrac{1}{10},\tfrac{2}{10},\tfrac{3}{10},\tfrac{4}{10},\tfrac{5}{10}\right) \nonumber\\
(c_{\beta_1},c_{\beta_2},c_{\beta_3},c_{\beta_4},c_{\beta_5})
=\mathstrut & (0.03150,0.18061,0.36563,0.31421,0.10801) .
\end{align}
Recall that all decimal numbers are truncations of the
actual values; one requires much higher precision than the displayed
numbers in order to obtain the answers below.
With the choices in~\eqref{eqn:5weights}, after substituting the known
Dirichlet coefficients
and then using the Ramanujan bound, we find
\begin{align}\label{eqn:errorexample}
Z&(s,\Upsilon_{20},\stan) = \nonumber \\
\mathstrut & 3.03930\,70864\,89527\,82778
+ 2.688\cdot 10^{-19} b_{83} + \cdots - 1.147\cdot 10^{-16} b_{107}+ \nonumber\\
& + \cdots -5.291 \cdot 10^{-18} b_{137} + \cdots + 1.216 \cdot 10^{-23} b_{199} + \cdots \nonumber\\
 =&\mathstrut  3.03930\,70864\,89527\,827 \pm 4.73 \cdot 10^{-15}.
\end{align}
Thus, by averaging only 5 evaluations of the $L$-function, the error
decreased by a factor of almost $10^{-5}$.

The weights $c_{\beta_j}$ in~\eqref{eqn:5weights} were determined
by finding the least-squares fit to
\begin{equation}\label{eqn:leastsquares}
\sum_{n:\ b_n\ \text{unknown}} Ram(b_n)^2
  \biggl(\sum_j c_{\beta_j} \delta_n(\beta_j,\tfrac12+ 10 i) \biggr)^2 =0,
\end{equation}
where $Ram(b_n)$ is the Ramanujan bound~\eqref{eqn:ram} for~$b_n$,
subject to $\sum c_{\beta_j} = 1$.
Note that in our actual examples the vast majority of unknown coefficients
have prime index, so the $ Ram(b_n)$ weighting is not important,
but we include it for completeness.
For the calculations in this paper, we use the
$n<1000$ for which $b_n$ is unknown in \eqref{eqn:leastsquares}.  

The error estimate in \eqref{eqn:errorexample} is an $L^1$ estimate,
not the $L^2$
estimate as shown in \eqref{eqn:leastsquares}, so actually it is possible to choose slightly
better weights than those used in our example.  
In Section~\ref{sec:lp} we show how to obtain the optimal result that
can arise from combining different evaluations
of the $L$-function, but for now we merely wish to illustrate the
seemingly surprising fact that appropriately combining several evaluations
can vastly decrease the error.

\subsection{Results}

\begin{figure}[htp] \begin{center}
\scalebox{0.7}[0.7]{\includegraphics{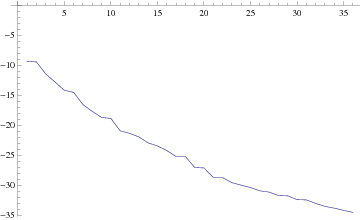}}
\caption{\sf 
The error obtained from a least-squares minimization
of the error for
combining
$n$ evaluations of 
$Z(\frac12+10 i,\Upsilon_{20},\stan)$ using the weight
functions $g(s)= e^{- i \beta s}$ with $\beta=j/10$.
The horizontal axis is $n$ and the vertical axis is
$\log_{10}$ of the resulting error estimate.
 \label{fig:stan_errors}}
\end{center}
\end{figure}

In Figure~\ref{fig:stan_errors} we plot the error obtained by combining
varying numbers of weight functions, where we evaluate
$Z(\frac12+10 i,\Upsilon_{20},\stan)$ with a weight function $g(s)=e^{i \beta s}$
with $\beta = j/10$ with $-10\le j \le 25$.
The horizontal axis is the number of terms averaged, where we start
with $\beta=\tfrac12$ and first use those $\beta$ which are closest to~$\frac12$.
The vertical axis is the error estimate on a $\log_{10}$ scale.
The lowest point on the graph, when we average all 36 evaluations,
corresponds to
\begin{equation}
\begin{split}Z(\tfrac12+&10 i,\Upsilon_{20},\stan)
=\\
&3.03930\,70864\,89528\,48108\,24603\,28442\,22509\,10
\pm
2.79 \times 10^{-35} .
\end{split}
\end{equation}

It is not clear from Figure~\ref{fig:stan_errors} whether one
would expect to obtain an arbitrarily small error by combining
sufficiently many test functions in the approximate functional
equations.  This is discussed in Section~\ref{sec:speculation}.

We briefly describe the corresponding calculations  for
$Z(\frac12+5 i,\Upsilon_{20},\adj)$.  Recall that, as shown
in Figure~\ref{fig:adj5error}, with a single test function of the
standard form we are not able to determine whether
that value is positive or negative.  Now we combine
five evaluations in the analogous way:
\begin{equation}
Z(s,\Upsilon_{20},\adj) = \sum_{j=1}^J c_{\beta_j} Z(s,\Upsilon_{20},\adj,\beta_j)
\end{equation}
where $\sum c_{\beta_j} =1$.  
Here the weight function is $g(s) = e^{-i \beta s + (s-5i)^2/500}$.

We make the specific choices
\begin{align}\label{eqn:10weights}
s=\mathstrut & \tfrac12+5 i\nonumber\\
(\beta_1,\beta_2,\beta_3,\beta_4,\beta_5)
   =\mathstrut & \left(\tfrac{3}{5},\tfrac{6}{5},\tfrac{9}{5},\tfrac{12}{5},\tfrac{15}{5}\right) \nonumber\\
(c_{\beta_1},c_{\beta_2},c_{\beta_3},c_{\beta_4},c_{\beta_5})
=\mathstrut & (0.035863,0.33504,0.47934,0.13827,0.01146) .
\end{align}
The result is
\begin{equation}
Z(\tfrac12 + 5 i,\Upsilon_{20},\adj) =
0.01556 \pm 0.0049.
\end{equation}
So we have determined that $Z(\tfrac12 + 5 i,\Upsilon_{20},\adj)$ is positive,
but we can only be certain of one significant figure in its decimal expansion.
Using this method, the best result we were able to obtain,
by averaging 11 evaluations,  is
\begin{equation}
Z(\tfrac12 + 5 i,\Upsilon_{20},\adj) =
0.01558768 \pm 0.00016.
\end{equation}
Averaging more weight function actually makes the result worse.
The explanation is simple:  since the adjoint $L$-function has high
degree, the error terms $\delta_n$ decrease very slowly.
The least-squares fit does not properly take into account
the contribution of a large number of small terms, so the least-squares
fit actually has a large $L^1$ norm. 

In the next section we describe some ``sanity checks'' on our method,
as suggested by the referee. Then, in Section~\ref{sec:lp}, we introduce
a different method which avoids some of the shortcomings
in the least-squares method.

\subsection{The method applied to known examples}\label{sec:realitychecks}

We check that our method gives correct results in some cases where
it is possible to evaluate the $L$-function using another method.

First we consider $L(s, \Delta)$, where $\Delta$
is the (unique) weight 12 cusp form for $SL(2,\Z)$,
which satisfies the functional equation
\begin{equation}
\Lambda(s,\Delta) := \Gamma_\C(s+\tfrac{11}{2}) L(s,\Delta)
= \Lambda(1-s,\Delta) .
\end{equation}
We will evaluate $Z(\tfrac12+100 i,\Delta)$
\emph{without using any Dirichlet coefficients}, other than
the leading coefficient~1.
That is, all we know about the $L$-function is its functional
equation and the fact that its Dirichlet coefficients
satisfy the Ramanujan bound.

In the approximate functional equation we will use the test
function $g_\beta(s)=e^{-i \beta s + (s-100i)^2/100}$, for
$\beta=-\frac{30}{20},-\frac{29}{20},...,\frac{69}{20},\frac{70}{20}$.
That is a total of 101 evaluations.  Using the least-squares method described
previously, we find 101 coefficients $c_\beta$ with $\sum c_\beta=1$.
Forming the weighted sum of the 101 evaluations of $Z(\tfrac12+100i,\Delta)$
and estimating the unknown terms as described previously,
we find
\begin{equation}
Z(\tfrac12+100i,\Delta) =
-0.23390\,65915\,56845\,20570\,65824\,17137\,27923\,81141\,00783
\pm 3.28 \times 10^{-42}.
\end{equation}
The given digits are correct to the claimed accuracy: 
the last three digits shown should be 880, and the actual
difference between the calculated value and the true value
is $9.66\times 10^{-44}$.

Next we consider the $L$-function associated to a weight 24
cusp form for $SL(2,Z)$.  Note that $S_{24}(SL(2,\Z)$ is two
dimensional, and every cusp form $f$ in that vector space
satisfies the functional equation
\begin{equation}\label{eqn:fe24}
\Lambda(s,f) := \Gamma_\C(s+\tfrac{23}{2}) L(s,f)
= \Lambda(1-s,f) .
\end{equation}
We will attempt to evaluate $Z(\frac12+100i,f)$ using
as few coefficients as possible.
Because there is more than one function satisfying~\eqref{eqn:fe24},
it seems obvious that we cannot evaluate such an $L$-function
without knowing any coefficients.

If we assume the cusp form $f$ is a Hecke eigenform, then
the Dirichlet coefficient $b_2$ determines the coefficients
$b_4$, $b_8$, etc, and it also allows us to eliminate every
even-index Dirichlet coefficient as an ``unknown.''
Since that does not seem like an adequately strenuous test of the method,
instead we will
assume nothing about the Dirichlet series except the
functional equation~\eqref{eqn:fe24} and a bound on the Dirichlet
coefficients.  We will assume a Ramanujan bound of the form
$|b_n| \le C_f d(n)$, where 
$d(n)$  is the divisor function and
$C_f$ is a constant depending only on the
cusp form~$f$. 
Is $f$ is a Hecke eigenform then $C_f=1$, and if
$f=A f_1+B f_2$ where $f_1,f_2\in S_{24}(SL(2,\Z)$ are the Hecke eigenforms,
then
$C_f = |A|+|B|$.

Using the same 101 test functions as in the case of $L(s,\Delta)$,
and choosing 101 weights to minimize the contribution of
$b_3$, $b_4$, ..., we find
\begin{align}\label{eqn:wt24at100i}
Z(\tfrac12 + 100i,f) =\mathstrut & 1.87042\,65340\,29268\,89914\,33391\,93910\,89610\,35060\,87410\ b_1 \cr
&+ 1.12500\,88863\,02338\,48447\,34844\,21487\,86375\,36206\,60254\ b_2 \cr
&\pm C_f \times 2.86 \times 10^{-43}.
\end{align}
Thus, we can evaluate $Z(\frac12+100i,f)$ to 42 decimal places,
knowing (up to a normalizing constant) only one Dirichlet coefficient.
The values of $b_2$ for the two Hecke eigenforms,
in the analytic normalization, are
\begin{equation}
b_2 = \frac{540 \pm 12 \sqrt{144169}}{2^{\frac{32}{2}}},
\end{equation}
and inserting those values finds that \eqref{eqn:wt24at100i} is correct.

Our third check on the method is to evaluate a degree 10 $L$-function
which can be also evaluated in an independent way.  To give a reasonable
match with the case of $L(\frac12+5 i, \Upsilon_{20},\adj)$,
we consider $L(\frac12+5 i,f)^5$, where $f$ is a weight 24 cusp
forms for $SL(2,\Z)$.  In other words, the same $L$-function as in
the previous example, except that we take its 5th power and evaluate
at $\frac12 + 5 i$.  
Note that this time there are two $L$-functions with the given
functional equation, with known values:
\begin{align}
L(\tfrac12+5 i,f_1)^5=\mathstrut & (-3.0527819)^5 = -265.14223\\
L(\tfrac12+5 i,f_2)^5=\mathstrut & (-0.7404879)^5 = -0.2226331 \ .
\end{align}

We will assume that we know the Euler factors
up through $p=79$, just as in the case of $L(\frac12+5 i, \Upsilon_{20},\adj)$.
If we only use a single test function of the standard form,
then the best error we can obtain is comparable to what we found in the
previous degree~10 case:
\begin{align}
L(\tfrac12+5 i,f_1)^5=\mathstrut & -265.204\pm 0.314 \\
L(\tfrac12+5 i,f_2)^5=\mathstrut & -0.22193 \pm 0.233\ .
\end{align}

Using the test functions $g_\beta(s) = e^{-i \beta s + (s-5i)^2/500}$,
selecting the 7 values of $\beta$ in the set
$\{-\frac{10}{10}, -\frac{2}{10},\frac{4}{10},\frac{6}{10},\frac{11}{10},
\frac{17}{10},\frac{18}{10},\frac{25}{10}\}$,
and using the least-squares method to find suitable weights, we find
\begin{align}
L(\tfrac12+5 i,f_1)^5=\mathstrut & -265.14224\pm 0.00117 \\
L(\tfrac12+5 i,f_2)^5=\mathstrut & -0.222664 \pm 0.00186\ .
\end{align}
Averaging only 7 evaluations decreased the error by a factor of more than 100,
and we see that the calculated values are in fact correct.
This confirms that our method gives consistent results in cases of
comparably complicated $L$-functions for which it is possible to give
an independent check on the calculations.

\section{Linear programming}\label{sec:lp}
In this section we view the evaluation of the $L$-function as an
optimization problem.  For example, we can view the equality of
the expressions in~\eqref{eqn:ex1} and~\eqref{eqn:ex2} as a
\emph{constraint}
on the value of the $L$-function.
Thus, the same calculations which were used as input for the
least-squares method described in the previous section can also be used
as input to a linear programming problem.

We set up the linear programming problem in the following way.
Let $Z(\frac12+10 i,\Upsilon_{20},\stan,g_j)$ denote the
evaluation of $Z(\frac12+10 i,\Upsilon_{20},\stan)$ using the weight
function $g_j$ in the approximate functional equation.
One evaluation, $Z(\frac12+10 i,\Upsilon_{20},\stan,g_1)$, is taken
as the objective.  The other evaluations are taken in pairs
and 
\begin{equation}\label{eqn:constraint}
Z(\tfrac12+10 i,\Upsilon_{20},\stan,g_1) -
Z(\tfrac12+10 i,\Upsilon_{20},\stan,g_j) = 0
\end{equation}
is interpreted as a constraint. The other constraints come
from the Ramanujan bound on the unknown coefficients.

In the above description there are infinitely many unknowns
and constraints. We eliminate the unknown coefficients $b_n$
with $n>1000$ by using the Ramanujan bound, replacing
the equality
\eqref{eqn:constraint} 
by a pair of inequalities.
Thus, we have a straightforward linear programming problem,
which we use to determine the minimum and maximum possible
values of $Z(\frac12+10 i,\Upsilon_{20},\stan)$.

\begin{figure}[htp]
\begin{center}
  \scalebox{0.70}[0.70]{\includegraphics{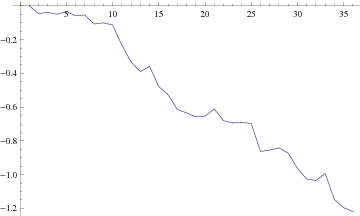}}
  \caption{\sf
The ratio of the errors in 
the linear programming and least-squares methods for
evaluating $Z(\frac12+10 i,\Upsilon_{20},\stan)$
using $n$ equations.
The horizontal axis is $n$ and the vertical axis is $\log_{10}$
of the ratio of errors.
\label{fig:LPstan10}
 }
\end{center}
\end{figure}

We implemented this idea using the same set of test functions
we used for the least-squares method.  The calculations were done
in Mathematica \cite{mathematica}, using the built-in
{\tt LinearProgramming} function with the {\tt Method -> Simplex}
option.  In all cases the linear programming method gave
better results, but not spectacularly better.
Figure~\ref{fig:LPstan10} shows the ratio of the errors
in the results of the two methods, on a $\log_{10}$ scale.
For example, when using 30 equations the error from the
linear programming approach was approximately $1/10$ the error
from the least-squares method.

\subsection{Computing unknown coefficients}

An interesting side-effect of the linear programming approach is that
it allows us to obtain information about the unknown coefficients.
Instead of treating the value of the $L$-function as the
objective, we can use an unknown coefficient 
as the objective.  Note that all the other constraints in the
problem are unchanged. Using the same method as described above
we find the following coefficients of the standard
$L$-function of~$\Upsilon_{20}$:
\begin{align}
b_{83}=\mathstrut& 0.48845\,58312\,724 \pm 2.4\times 10^{-12}\cr
b_{89}=\mathstrut& 0.10561\,760640 \pm 2.7\times 10^{-10}\cr
b_{97}=\mathstrut& 0.46813\,5808 \pm 1.5 \times 10^{-7}
\end{align}
Since the eigenvalues of $\Upsilon_{20}$ are integral and of a known
size, if we were to know $b_{83}$ to 35 digits, we would determine it
exactly.  This is computationally expensive but perhaps not as
expensive as if we were to compute more Fourier coefficients of $\Upsilon_{20}$.

These results can be checked once more Fourier coefficients of
$\Upsilon_{20}$ are computed.  By the $n^4$ argument given in in the
introduction,
computing exact values of 
$b_p$ for $p\le 97$ will take about twice as much work as it took
to compute those for $p\le 79$.

\section{Conclusions and further questions}\label{sec:speculation}

We have shown that, at the cost of a lot of computation, one can evaluate
an $L$-function to high precision using only a small number of
coefficients.  That this is theoretically possible is not surprising:
$L$-functions are very special objects, and the data we have
for the $L$-function considered here (the functional equation and the
first several coefficients) presumably specify the $L$-function uniquely.
Thus, in an abstract sense there is no new information in the missing
coefficients. But the question remains as to whether our methods
accomplish this in practice.

\begin{question} Can the method of calculating an $L$-function
by evaluating the approximate functional equation and then
averaging to minimize the contributions of the unknown
coefficients, determine numerical values of the $L$-function
to arbitrary accuracy?
\end{question}

Because the approximate functional equation requires a huge
number of terms to evaluate a high-degree $L$-function, it would
be significant if the weights we obtained by our methods 
could be determined without
actually calculating all the terms with unknown coefficients.

\begin{problem} Devise a method of determining an optimal weight
function in the approximate functional equation without first
calculating a large number of terms which do not actually
contribute to the final answer.
\end{problem}

Implicit in the above problem is the requirement that 
one knows the calculated value as well as an estimate of the
error.

\begin{problem} Is there any meaning to the weights determined
by the least-squares method?
\end{problem}

The weights which appear in our least-squares method depend on the
point at which the $L$-function is evaluated.  It might be helpful
to consider a case where there are two $L$-functions with the same
functional equation, such as the spin $L$-functions of
$\Sp(4,\Z)$ Siegel modular forms of weight $k\ge 22$.

Without progress on these problems, or a completely new method,
there is little hope of making extensive computations of
high-degree $L$-functions.

\bibliography{evaluation_few3}
\bibliographystyle{plain}

\affiliationone{David W. Farmer\\
American Institute of Mathematics\\
360 Portage Ave\\
Palo Alto, CA 94306\\
USA\\
\email{farmer@aimath.org}
}
\affiliationtwo{Nathan C. Ryan\\
Department of Mathematics\\
Bucknell University\\
Lewisburg, PA 17837\\
USA\\
\email{nathan.ryan@bucknell.edu}
}

\end{document}